\documentclass{amsart}
\usepackage{amsfonts}
\usepackage{amsmath,amscd}
\usepackage{amssymb}
\usepackage{amsthm}
\usepackage{newlfont}
\newcommand{\f}{\frac}

\newcommand{\ds}{\displaystyle}

\begin{document}





\title[Relative $n$-isoclinism Classes and Relative $n$-th Nilpotency Degree of
Finite Groups]
 {Relative $n$-isoclinism Classes and Relative $n$-th Nilpotency Degree of
Finite Groups}

\author[A. Erfanian]{Ahmad Erfanian}
\address{Department of Mathematics, Faculty of Mathematical Sciences\\
Ferdowsi University of Mashhad, Mashhad, Iran.}
\email{erfanian@math.um.ac.ir}

\author[R. Rezaei]{Rashid Rezaei}

\address{Department of Mathematics, Faculty of Mathematical Sciences\\
University of Malayer\\
Post Box: 657719--95863, Malayer, Iran} \email{ras$\_$rezaei@yahoo.com}

\author[F.G. Russo]{Francesco G. Russo}
\address{Department of Mathematics\\
University of Naples Federico II\\
via Cinthia, 80126,  Naples, Italy}
\email{francesco.russo@dma.unina.it}



\date{\today}
\dedicatory{}

\begin{abstract}
The purpose of the present paper is to consider the notion of isoclinism between two finite groups and its
generalization to $n$-isoclinism, introduced by J. C. Bioch in 1976. A weaker form of $n$-isoclinism, called
relative $n$-isoclinism, will be discussed. This will allow us to improve some classical results in literature.
We will point out the connections between a relative $n$-isoclinism and the notions of commutativity degree,
$n$-th nilpotency degree and relative $n$-th nilpotency degree, which arouse interest in the classification of
groups of prime power order in the last years.
\end{abstract}

\maketitle















\section{Introduction}
\indent  The notion of isoclinism was introduced by P. Hall in 1940
in [8], looking for a satisfactory classification of groups of prime
power order. However, the notion of isoclinism holds both in finite
and infinite group theory. It is obvious that an isomorphism is an
equivalence relation in the class of all groups and this allows us
to classify two groups. An isoclinism is a more general equivalence
relation in the class of all groups and it is easy to see that two
abelian groups fall into the same equivalence class with respect to
isoclinisms (see [1, Theorem 1.4] or [8] or [9]). Roughly speaking,
two groups $H$ and $K$ are isoclinic if their central quotients
$H/Z(H)$, $K/Z(K)$  are isomorphic and if their commutator subgroups
$H'$, $K'$ are isomorphic.  If we look at the construction of the
finite extra-special $2$-groups (see [15, pp.145--147]) and at the
construction of the quaternion groups (see [15, pp.140--141]), then
we will find such groups in the situation which has been just
described. For instance, we may think to the dihedral group $D_8$ of
order 8 and to the quaternion group $Q_8$ of order 8. We note that
both $D_8/Z(D_8)$, $Q_8/Z(Q_8)$ are isomorphic and $D'_8$, $Q'_8$
are isomorphic. \

Situations as we just mentioned have been largely studied in
literature and it is well known the role of the isoclinism in the
classification of the finite groups of prime power order as we can
see  in [3, 7, 8, 9, 10, 11, 12, 13, 14].

Such references and many other works of the same authors investigate
those group theoretical properties which are invariant under
isoclinism. For example, J. C. Bioch and R. W. van der Waall [2]
proved  the invariance under isoclinism of the following hierarchy
of classes of finite groups:\\
abelian $<$ nilpotent $<$ supersoluble  $<$ strongly-monomial $<$
monomial $<$ soluble.\

A successive contribution in the classification of the groups of
prime power order with respect to the notion of isoclinism was given
already by P. Hall. He introduced the varieties of groups in [9]. We
can easily see that the variety of all trivial groups and the
variety of all abelian groups yield to the notion of isomorphism and
isoclinism, respectively (see [1, 2, 10]). Varieties which extend
the variety of abelian groups yield to notions which extend that of
isoclinism . In this context, J. C. Bioch took the variety of all
nilpotent groups of class at most $n$ ($n$ is a positive integer)
and introduced the notion of $n$-isoclinism of groups (see [1]).
Successive contributions to the works of J. C. Bioch were given by
[10].

Recently, some interesting connections have been found between the
notion of isoclinism and the probability that two randomly chosen
elements of a group permute. In particular, probabilistic bounds
have been found in [3, 4, 5, 6, 12, 13, 16], giving restrictions in
terms of the structure of a finite group. This paper follows such a
line of research. Our main results extend both some classical
results in literature on isoclinic groups as [1, 2, 8, 9, 10, 11]
and some of the results in [3, 4, 5, 6, 12, 13, 16]. Section 2
recalls the technical definitions which we will use in Section 3 for
getting to our main results. Some open questions both in infinite
and finite case have been shown at the end of Section 3. Our
notation and terminology is standard and referred to [15].

\section{Preliminaries}
\indent In this section, we recall two important definitions.
Firstly, we recall the notion of commutativity degree and its
generalizations. Secondly, we recall the notion of isoclinism, and
its generalizations. They can be found in [1, 2, 3, 4, 8, 9, 10,
13]. From this point, the symbol $n$ will denote always a positive
integer.

For any finite group $G$, \textit{the commutativity degree of }$G$
is defined as the number of pairs $(x,y)$ in $G\times G$ such that
$xy=yx$, divided by $|G|^2$.  We denote it by $d(G)$. In symbols:
$$d(G)=\ds\f{|\{(x,y)\in G^2 : [x,y]=1\}|}{|G|^2}.
\hspace{5.2cm}(2.1)$$ It is clear that $d(G)=1$ if and only if $G$
is abelian; and for any finite non-abelian group $G$, $d(G)\leq\ds\f
 58$. Furthermore, this bound is achieved
if and only if $G/Z(G)$ is a $2$-elementary abelian group of rank
$2$. Further details can be found in [4, Theorem 2] and [6,
Introduction].\

A way of extending $d(G)$ to $(n+1)$-tuples $(x_1,x_2,...,x_{n+1})$
in $G^{n+1}$ is to consider the number of such $(n+1)$-tuples with
the property that $[x_1,x_2,...,x_{n+1}]=1$, divided by the order of
$G^{n+1}$. This is denoted by $d^{(n)}(G)$ and called the
\textit{$n$-th nilpotency degree of }$G$. In
symbols:$$d^{(n)}(G)=\ds\f{|\{(x_1,x_2,...,x_{n+1})\in G^{n+1} :
[x_1,x_2,...,x_{n+1}]=1\}|}{|G|^{n+1}}. \hspace{1cm}(2.2)$$ It is
obvious that $d^{(n)}(G)$ extends $d(G)$ and $d^{(1)}(G)=d(G)$ (see
[13, Section 4] and [3] for more details).

Now, we will introduce two further generalizations of $d(G)$. If we
take a subgroup $H$ of a finite group $G$, we may define the
\textit{relative commutativity degree of $H$ in $G$}, denoted by
$d(H,G)$, as $$d(H,G)=\ds\f{|\{(h,g)\in H\times G :
[h,g]=1\}|}{|H||G|}\hspace{4.6cm}(2.3)$$ and we may define the
 \textit{relative n-th nilpotency degree of $H$ in  $G$}, denoted by
$d^{(n)}(H,G)$, as $$d^{(n)}(H,G)=\ds\f{|\{(x_1,x_2,...,x_n,g)\in
H^n\times G : [x_1,x_2,...,x_n,g]=1 \}|}{|H|^n|G|}.
\hspace{0.5cm}(2.4)$$ If $H=G$, then $d(H,G)=d(G)$ and
$d^{(n)}(H,G)=d^{(n)}(G)$ as we see in [4, Section 3].

The second part of this section is devoted to define isoclinism,
$n$-isoclinism, relative isoclinism and relative $n$-isoclinism. We
start with defining isoclinism between two groups $H$ and $G$,
following [8] and recent considerations in [13, Section 2]. \\

{\underline{\bf{Definition 2.1.}}}~ \textit{Let $G$ and $H$ be two
groups; a pair $(\varphi,\psi)$ is called an
isoclinism from $G$ to $H$ if :\\
  (i) $\varphi$ is an isomorphism from $\ds \f {G}{Z(G)} $ to $\ds \f
  {H}{Z(H)};$\\
  (ii) $\psi$ is an isomorphism from $G'$ to $H'$ ;\\
  (iii) the following diagram is commutative:\\
\[\begin{CD}
 \ds \f G{Z(G)}\times \f {G}{Z(G)}@>\varphi>> \ds\f H{Z(H)}\times \f
H{Z(H)} \\
@Va_GVV @Va_HVV\\
 G'@>\psi>> H',
\end{CD}\]}
\textit{where $a_G( g_1 Z(G),g_2 Z(G))=[g_1,g_2]$ and $a_H( h_1
Z(H),h_2 Z(H))=[h_1,h_2]$, for each $g_1,g_2 \in G$ and $h_1,h_2 \in
H$.}\\
\\
If there is an isoclinism from $G$ to $H$, we say that $G$ and $H$
are $isoclinic$, writing briefly $G\sim H$. It can be easily checked
that the relation $\sim$ given in Definition 2.1 is an equivalence
relation. Moreover, it is obvious that if $G$ and $H$ are
isomorphic, then they are isoclinic. But the converse is not true.
By a simple calculation, one may easily see that $Q_8$ and $D_8$
are isoclinic but not isomorphic.\\
A simple relation between isoclinism and commutativity degree is
that two isoclinic finite groups have the same commutativity degree
as we see in [13, Lemma 2.4]. Another relation between finite groups
of commutativity degree equal to $\frac{1}{2}$ and groups which are
isoclinic to the symmetric group $S_3$ has been given in [13,
Theorem 3.1].

Now, we give the notion of $n$-isoclinism between two groups,
generalizing Definition 2.1. This notion has been extensively
investigated in [10, Sections 3, 5, 7]. \\

\indent {\underline{\bf{Definition 2.2.}}}~ \textit{Let $G$ and $H$
be two groups; a pair $(\alpha,\beta)$ is called  n-isoclinism from
$G$ to $H$
if :\\
 (i) $\alpha$ is an isomorphism from $\ds \f G{Z_n(G)} $ to $\ds \f
 H{Z_n(H)}$;\\
 (ii) $\beta$ is an isomorphism from $\gamma_{n+1}(G)$ to
$\gamma_{n+1}(H)$; \\
 (iii) the following diagram is commutative:\\
\[\begin{CD}
\ds \f G{Z_n(G)}\times ...\times\f G{Z_n(G)}
@>\alpha^{n+1}>> \ds \f H{Z_n(H)}\times...\times\f H{Z_n(H)}\vspace{.3cm}\\
@V\gamma(n,G)VV @V\gamma(n,H)VV\vspace{.3cm}\\
\gamma_{n+1}(G)@>\beta>>\gamma_{n+1}(H),
\end{CD}\]
where $$\gamma (n,G)((g_1Z_n(G),...,g_nZ_n(G),g_{n+1}Z_n(G)))=[g_1,
..., g_n,g_{n+1}]$$ and $$\gamma
(n,H)((h_1Z_n(H),...,h_nZ_n(H),h_{n+1}Z_n(H)))=[h_1, ...,
h_n,h_{n+1}],$$ for each $g_1,...,g_n,g_{n+1} \in G$ and
$h_1,...,h_n,h_{n+1} \in H$.}\\
\\
One may easily check that the maps $\gamma(n,G)$ and $\gamma(n,H)$
in Definition 2.2 are well-defined.

If there is an $n$-isoclinism between $G$ and $H$, we say that $G$ and $H$ are $n$-$isoclinic$, writing briefly
by $G \ _{\widetilde{n}} \ H$. It is clear that $\ _{\widetilde{_1}}$\ and $\sim$ coincide, that is, a
$1$-isoclinism is an isoclinism. Follows from Definition 2.2 that $_{\widetilde{n}}$ is an equivalence relation
(see [1] or [10, Section 3]). If $G$ and $H$ are $n$-isoclinic, then they are $(n+1)$-isoclinic by [10, Theorem
5.2]. This will be extended in Section 3.

Note that our terminology can be found in [1, 2, 7, 8, 9, 11], as
most of the terminology in the present section.

We know that Definitions 2.1, 2.2 and  Equations (2.1), (2.2) are
correlated by results as [13, Lemma 2.4] or [3, Theorem B]. Indeed,
[3, Theorem B] states that two $n$-isoclinic groups have the same
$n$-th nilpotency degree.

Now, we state the last definition of the present section which
recalls [10, Definition 7.1] and extends Definitions 2.1 and 2.2.\

\vspace{0.4cm} \indent {\underline{\bf{Definition 2.3.}}}~
\textit{Let $G_1$, $G_2$ be two  groups, $H_1$  a subgroup of $G_1$
and $H_2$  a subgroup of $G_2$. A pair $(\alpha,\beta)$ is said
to be a relative $n$-isoclinism from $(H_1,G_1)$ to $(H_2,G_2)$ if we have the following conditions:\\
(i) $\alpha$ is an isomorphism from $G_1/Z_n(G_1)$ to $G_2/Z_n(G_2)$
such that \ the restriction of $\alpha$ under $H_1/(Z_n(G_1)\cap
H_1)$ is an isomorphism from $H_1/(Z_n(G_1)\cap H_1)$ to
$H_2/(Z_n(G_2)\cap
H_2)$;\\
  (ii) $\beta$ is an isomorphism from $[_nH_1,G_1]$ to
  $[_nH_2,G_2]$;\\
  (iii) the following diagram is commutative:\\ {\small
\[\begin{CD}
 \f{H_1}{Z_n(G_1)\cap H_1}\times ...\times\f{H_1}{Z_n(G_1)\cap
H_1}\times \f {G_1}{Z_n(G_1)}
  @>\alpha^{n+1}>> \f {H_2}{Z_n(G_2)\cap H_2}\times...\times\f
{H_2}{Z_n(G_2)\cap H_2}\times \f {G_2}{Z_n(G_2)} \vspace{.3cm}\\
@V\gamma(n,H_1,G_1) VV @V\gamma(n,H_2,G_2)VV \vspace{.3cm}\\
[_nH_1,G_1] @>\beta>>[_nH_2,G_2].
\end{CD}\]
where
$$\gamma (n,H_1,G_1)((h_1(Z_n(G_1)\cap H_1),...,h_n(Z_n(G_1)\cap H_1),g_1Z_n(G_1)))=[h_1,
..., h_n,g_1]$$ and $$\gamma (n,H_2,G_2)((k_1(Z_n(G_2)\cap
H_2),...,k_n(Z_n(G_2)\cap H_2),g_2Z_n(G_2)))=[k_1, ..., k_n,g_2],$$
for each $h_1,...,h_n \in H_1$, $k_1,...,k_n \in H_2$, $g_1 \in
G_1$, $g_2\in G_2$.}}\\
\\
It is easy to check that  the maps $\gamma(n,H_1,G_1)$ and
$\gamma(n,H_2,G_2)$ in Definition 2.2 are well-defined.

If Definition 2.3 is satisfied, we say that $(H_1,G_1)$ and $(H_2,G_2)$ are $relative$ $n$-$isoclinic$, writing
briefly $\ds(H_1,G_1) \ _{\widetilde{n}} \ (H_2,G_2).$ Follows from Definition 2.3 that $_{\widetilde{n}}$ is an
equivalence relation (see also [1] or [10, Section 3]). Note that there is not ambiguity if we use the symbol
$_{\widetilde{n}}$ both in Definitions 2.2 and 2.3, because in Definition 2.2 $_{\widetilde{n}}$  is referred to
groups and in Definition 2.3 $_{\widetilde{n}}$  is referred to pairs of groups. We should also note that if two
pairs $(H_1,G_1)$ and $(H_2,G_2)$ are relative $n$-isoclinic, then it is not necessary that $G_1$ and $G_2$ are
$n$-isoclinic. For instance, assume that $SL(2,5)$ is the special linear group of order 120 and $PSL(2,5)$ is
the projective special linear group of order 60. Then we can observe that $(Z(SL(2,5)), SL(2,5))$ and $(1,
PSL(2,5))$ are relative 1-isoclinic. On another hand, we note that $|[SL(2,5),SL(2,5)]|=120$ and
$|[PSL(2,5),PSL(2,5)]|=60$ so the relative 1-isoclinism is not an isoclinism.

Now, we end the present section, stating our main results.

\vspace{0.4cm} {\underline{\bf{Theorem A.}}}~ \textit{Let $G_1$ and
$G_2$ be two $n$-isoclinic groups. For every subgroup $H_1$ of
$G_1$, there exists a subgroup $H_2$ of $G_2$ such that $H_1$ and
$H_2$ are $n$-isoclinic.}\\

\vspace{0.4cm} {\underline{\bf{Theorem B.}}}~ {\em Let $G_1$ and
$G_2$ be two finite groups, $H_1$ be a subgroup of $G_1$ and $H_2$
be a subgroup of $G_2$. If $\ds(H_1,G_1) \ _{\widetilde{n}} \
(H_2,G_2)$, then $d^{(n)}(H_1,G_1)
=d^{(n)}(H_2,G_2)$}.\\

\vspace{0.4cm}{\underline{\bf{Theorem C.}}}~ \textit{Let $G$ be a
group and $H,N$ be subgroups of $G$ such that $N\triangleleft G$ and
$N\subseteq H$. Then for all $n\geq0$,
$$(\f HN,\f GN) \ _{\widetilde{n}} \ (\f{H}{N\cap\gamma_{n+1}(G)},\f{G}{N\cap\gamma_{n+1}(G)}).$$ In particular, if
$N\cap \gamma_{n+1}(G)=1$, then $(H,G) \ _{\widetilde{n}} \
(H/N,G/N)$.}\\

\vspace{0.4cm}{\underline{\bf{Theorem D.}}}~{\em Let $H$ be a
subgroup of a finite group $G$. Then the following
statements are true.\\
{\textit(i)} If $G=HZ_n(G)$, then
$d^{(n)}(H)=d^{(n)}(H,G)=d^{(n)}(G)$. \\
{\textit (ii)} $d^{(n)}(H,G)=d^{(n)}(\varphi(H),G)$ for every
$\varphi\in Aut(G)$.}\\

\vspace{0.4cm}{\underline{\bf{Theorem E.}}}~{\em Let $H$ be a
subgroup of a finite group $G$ such that $Z(G)\subseteq H$. Then
$d(H,G)=\f34$ if and only if $(H,G)$ and $( \langle a\rangle,D_8)$
are relative 1-isoclinic, where $\langle a\rangle$ is a subgroup of
order $4$ in $D_8$.}

\section{Proof of Main Theorems}\
\indent This section is devoted to describe our main results.\\

{\underline{\bf{Lemma 3.1.}}}~ \textit{ $\ds(H_1,G_1) \
_{\widetilde{n}} \ (H_2,G_2)$ if and only if there exist two
isomorphisms $\alpha$ and $\beta$ such that
$$\alpha:G_1/Z_n(G_1)\rightarrow G_2/Z_n(G_2),$$
$$\beta:[_nH_1,G_1]\rightarrow[_nH_2,G_2],$$ $\alpha(H_1/(Z_n(G_1)\cap
H_1))=H_2/(Z_n(G_2)\cap H_2) \ and \ for\ all \ g_1\in G_1 \ and \
h_i\in H_1,$
$$\beta([h_1,...,h_n,g_1])=[k_1,...,k_n,g_2],$$ where
$g_2\in\alpha(g_1Z_n(G_1))$, $k_i\in\alpha(h_i(Z_n(G_1)\cap H_1))$
and $1\leq i\leq n$.}\\
\\
\indent {\it \bf Proof.}~ It is clear by Definition 2.3. $\square$\\

\vspace{0.4cm} {\underline{\bf{Lemma 3.2.}}}~ \textit{Let $G_1$ and
$G_2$ be two n-isoclinic groups. Then for every subgroup $H_1$ of
$G_1$ there exists a subgroup $H_2$ of $G_2$ such that $(H_1,G_1)
\ _{\widetilde{n}} \ (H_2,G_2).$}\\
\\
\indent {\it \bf Proof.}~ Assume that $(\alpha,\beta)$ is an
$n$-isoclinism from $G_1$ to $G_2$ as in Definition 2.3. Put
\[\begin{array}{lcl}
H_2= \{x\in G_2 \ | \ \exists \ h\in H_1 \ \ \textrm{s.t.} \ \
\alpha(hZ_n(G_1))= xZ_n(G_2)\},
\end{array}\]
where $H_1$ is an arbitrary subgroup of $G_1$. It is clear that
$H_2$ is a subgroup of $G_2$ and $\alpha(H_1/(Z_n(G_1)\cap
H_1))=H_2/(Z_n(G_2)\cap H_2).$ If $\beta'$ denotes the restriction
of $\beta$ under $[_nH_1,G_1]$, then $\beta'$ is an isomorphism from
$[_nH_1,G_1]$ to $[_nH_2,G_2]$ and for all $g_1\in G_1$, $h_i\in
H_1$, $1\leq i\leq n$ we have
\[\begin{array}{lcl}
\beta'([h_1,...,h_n, g_1])=\beta([h_1,...,h_n, g_1])\\
=\beta\gamma(n, G_1)(h_1Z_n(G_1),...,h_nZ_n(G_1),g_1Z_n(G_1))\\
=\gamma(n, G_2)\alpha^{n+1}(h_1Z_n(G_1),...,h_nZ_n(G_1),g_1Z_n(G_1))\\
=[k_1,...,k_n,g_2],
\end{array}\]
where $g_2\in\alpha(g_1Z_n(G_1)\cap H_1)$ and
$k_i\in\alpha(h_iZ_n(G_1)\cap H_1)$. Therefore $(\alpha, \beta')$ is
a relative $n$-isoclinism from $(H_1,G_1)$ to $(H_2,G_2)$, by Lemma 3.1. $\square$\\

Now, we may prove Theorem A.\\

{\textit{\bf{Proof of Theorem A.}}}~ Assume that $G_1$ and $G_2$ are
two $n$-isoclinic groups  and $H_1$ is an arbitrary subgroup of
$G_1$. By Lemma 3.2, there exists a subgroup $H_2$ of $G_2$ such
that $(H_1,G_1) \ _{\widetilde{n}} \ (H_2,G_2).$ Let
$(\alpha,\beta)$ be a relative $n$-isoclinism from $(H_1,G_1)$ to
$(H_2,G_2)$ as in Definition 2.3. We have the isomorphisms
$$\alpha': H_1/(Z_n(G_1)\cap H_1)\rightarrow H_2/(Z_n(G_2)\cap H_2),$$
and
$$\beta':\gamma_{n+1}(H_1) \rightarrow \gamma_{n+1}(H_2),$$
where $\alpha'$ is the restriction of $\alpha$ under
$H_1/(Z_n(G_1)\cap H_1)$  and $\beta'$ is the restriction of $\beta$
under $\gamma_{n+1}(H_1)$. By [10, Lemma 3.11], if $A=Z_n(G_1) \cap
H_1\subseteq Z_n(H_1)$ and $B=Z_n(G_2) \cap H_2\subseteq Z_n(H_2)$,
then $H_1$ and $H_2$ are $n$-isoclinic. $\square$\\

The following fact compares two relative $n$-isoclinisms.\\

{\underline{\bf{Proposition 3.3.}}}~ \textit{If
$(H_1,G_1) \ _{\widetilde{n}} \ (H_2,G_2) $, then $(H_1,G_1) \ _{\widetilde{_{n+1}}} \ (H_2,G_2) $.}\\

\indent {\it \bf Proof.}~ Let $(\alpha,\beta)$ be a relative
$n$-isoclinism from $(H_1,G_1)$ to $(H_2,G_2)$. Define
\[ \varphi:\ds\f {G_1}{Z_{n+1}(G_1)}\rightarrow \f
{G_2}{Z_{n+1}(G_2)} \ \ \ \textrm{and}  \ \
\psi:[_{n+1}H_1,G_1]\rightarrow [_{n+1}H_2,G_2]
\] by the rules
$ \varphi(g_1Z_{n+1}(G_1))=g_2Z_{n+1}(G_2)$ and
$\psi([x,g_1])=[\beta(x),g_2]$,
 where
$\alpha(g_1Z_n(G_1))=g_2Z_n(G_2)$ and $x\in \gamma_{n+1}(H).$ Now,
we claim that $\varphi$ and $\psi$ are isomorphisms.

Assume that $aZ_{n+1}(G_1)=bZ_{n+1}(G_1)$ for some $a,b\in G_1$,
then $\varphi(aZ_{n+1}(G_1))=xZ_{n+1}(G_2)$ ,
$\varphi(bZ_{n+1}(G_1))=yZ_{n+1}(G_2)$, where $x,y\in G_2$.
Therefore either
\[\begin{array}{lcl} b^{-1}aZ_n(G_1)\in
\ds\f{Z_{n+1}(G_1)}{Z_n(G_1)}=Z(\f {G_1}{Z_n(G_1)})
\end{array}\] or \[\begin{array}{lcl} \alpha(b^{-1}aZ_n(G_1))\in \alpha(Z(\ds\f
{G_1}{Z_n(G_1)}))=Z(\ds\f
{G_2}{Z_n(G_2)})=\ds\f{Z_{n+1}(G_2)}{Z_n(G_2)}
\end{array}\] or
\[\begin{array}{lcl} y^{-1}xZ_n(G_2) \in
\ds\f{Z_{n+1}(G_2)}{Z_n(G_2)}
\end{array}\] i.e. $\varphi(aZ_{n+1}(G_1))=\varphi(bZ_{n+1}(G_1))$. Thus $\varphi$ is
well-defined.

If $\varphi(aZ_{n+1}(G_1))=\bar{1}$, then
$\varphi(aZ_{n+1}(G_1))=xZ_{n+1}(G_2)=\bar{1}$ for some $x\in G_2$.
Thus $x\in Z_{n+1}(G_2)$ and so,
\[\begin{array}{lcl} \alpha(aZ_n(G_1))&=&xZ_n(G_2)\in
\ds\f{Z_{n+1}(G_2)}{Z_n(G_2)}=Z(\ds\f {G_2}{Z_n(G_2)})\vspace{0.3cm}\\
&=&\ds\alpha( Z(\f
{G_1}{Z_n(G_1)}))=\alpha(\f{Z_{n+1}(G_1)}{Z_n(G_1)}).
\end{array}\] Hence, $a\in Z_{n+1}(G_1)$ and therefore $\varphi$ is
injective.\

From $\alpha$ surjective, we conclude that $\varphi$ is surjective.
Therefore $\varphi$ is an isomorphism and
$\varphi(H_1/(Z_{n+1}(G_1)\cap H_1))=H_2/Z_{n+1}(G_2)\cap H_2$.

Since $\beta$ is an isomorphism, $\psi$ is a monomorphism. $\psi$ is
surjective because, if $[k_1,...,k_{n+1} ,g_2]\in [_{n+1}H_2,G_2]$
then there exist elements $h_1, ...,h_{n+1} \in H_1$ and $g_1\in
G_1$ such that $\beta([h_1,...,h_{n+1}])=[k_1,...,k_{n+1}]$ and
$\alpha(g_1Z_n(G_1))=g_2Z_n(G_2)$. Therefore
$\psi([h_1,...,h_{n+1},g_1])=[\beta(h_1,...,h_{n+1}),g_2]=[k_1,...,k_{n+1},g_2]$,
and so  $\psi$ is an isomorphism.

Now Lemma 3.1 implies that $(\varphi,\psi)$ is a relative $n+1$-isoclinism from $(H_1,G_1)$ to $(H_2,G_2).\square$\\

{\textit{\bf{Proof of Theorem B.}}}~ Let $(\alpha,\beta)$ be a
relative $n$-isoclinism from $(H_1,G_1)$ to $(H_2,G_2)$. We have
\vspace{.3cm}\\
$|\ds\f {H_1}{Z_n(G_1)\cap H_1}|^n|\f {G_1}{Z_n(G_1)}| \
d^{(n)}(H_1,G_1) \vspace{.3cm} \\
 = \f 1{|Z_n(G_1)\cap H_1|^n|Z_n(G_1)|}|  \{ (h_1,...,h_n,g_1)\in
H_1^n\times G_1 : \ [h_1,...,h_n,g_1]=1 \}| \vspace{.3cm}\\
 =\f 1{|Z_n(G_1)\cap H_1|^n|Z_n(G_1)|}| \cdot \vspace{.3cm}\\ \{ (h_1,...,h_n,g_1)\in
H_1^n\times G_1 \ : \ \gamma(n,H_1,G_1)(h_1(Z_n(G_1)\cap
H_1),...,h_n(Z_n(G_1)\cap H_1),g_1Z_n(G_1))=1 \}| \vspace{.3cm} \\
= |\{(\bar{h_1},...,\bar{h_n},\bar{g_1}) \in \f {H_1}{Z_n(G_1)\cap H_1} \times...\times \f {H_1}{Z_n(G_1)\cap
H_1} \times \f {G_1}{Z_n(G_1)}  \ : \ \gamma(n,H_1,G_1)(\bar{h_1},...,\bar{h_n},\bar{g_1})=1 \}| \vspace{.3cm} \\
= |\{(\bar{h_1},...,\bar{h_n},\bar{g_1}) \in \f {H_1}{Z_n(G_1)\cap H_1} \times...\times \f {H_1}{Z_n(G_1)\cap
H_1} \times \f {G_1}{Z_n(G_1)}  \ : \ \beta(\gamma(n,H_1,G_1)(\bar{h_1},...,\bar{h_n},\bar{g_1}))=1 \}|$\vspace{.3cm}\\
(because $\beta$ is an isomorphism) \vspace{.3cm} \\
$= |\{(\bar{h_1},...,\bar{h_n},\bar{g_1}) \in \ds\f {H_1}{Z_n(G_1)\cap H_1} \times...\times \f
{H_1}{Z_n(G_1)\cap H_1} \times \f {G_1}{Z_n(G_1)} \ : \
\gamma(n,H_2,G_2)(\alpha(\bar{h_1}),...,\alpha(\bar{h_n}),\alpha(\bar{g_1}))=1
\}| \vspace{.3cm}$\\
(by commutativity of diagram as in Definition 2.3) \vspace{.3cm} \\
$= |\{(\bar{k_1},...,\bar{k_n},\bar{g_2}) \in \ds\f {H_2}{Z_n(G_2)\cap H_2} \times...\times\f {H_2}{Z_n(G_2)\cap
H_2}\times\f {G_2}{Z_n(G_2)} \ : \ \gamma(n,H_2,G_2)(\bar{k_1},...,\bar{k_n},\bar{g_2})=1 \}|$\vspace{.3cm}\\
(because $\alpha$ is an isomorphism) \vspace{.3cm} \\
$=\ds\f 1{|Z_n(G_2)\cap H_2|^n|Z_n(G_2)|}| \{ (k_1,...,k_n,g_2)\in
H_2^n\times G_2 \ : \ \gamma(n,H_2,G_2)(\bar{k_1},...,\bar{k_n},\bar{g_2})=1 \}| \vspace{.3cm}\\
= \f 1{|Z_n(G_2)\cap H_2|^n|Z_n(G_2)|}| \{ (k_1,...,k_n,g_2)\in
H_2^n\times G_2 \ : \ [k_1,...,k_n,g_2]=1 \}| \vspace{.3cm}\\
=|\f {H_2}{Z_n(G_2)\cap H_2}|^n|\f {G_2}{Z_n(G_2)}| \
d^{(n)}(H_2,G_2)
.$ \vspace{.3cm}\\
Therefore
\[\begin{array}{lcl} |\ds\f {H_1}{Z_n(G_1)\cap H_1}|^n|\f
{G_1}{Z_n(G_1)}| \ d^{(n)}(H_1,G_1) =|\f {H_2}{Z_n(G_2)\cap
H_2}|^n|\f {G_2}{Z_n(G_2)}| \ d^{(n)}(H_2,G_2)
\end{array}.\]

Since $(\alpha,\beta)$ is a relative $n$-isoclinism,
$d^{(n)}(H_1,G_1)=d^{(n)}(H_2,G_2)$ and the result follows.$\square$\\

Theorem C generalizes [1, Lemma 1.3].\\

{\textit{\bf{Proof of Theorem C.}}}~ Put $\bar{G}=G/N$ and
$\widetilde{G}=G/(N\cap[_nH,G])$. Since $\bar{g}\in Z_n(\bar{G})$ if
and only if $\widetilde{g}\in Z_n(\widetilde{G})$, the map $\alpha$
from $\bar{G}/Z_n(\bar{G})$ onto $\widetilde{G}/Z_n(\widetilde{G})$
given by
$\alpha(\bar{g}Z_n(\bar{G}))=\widetilde{g}Z_n(\widetilde{G})$ is an
isomorphism and
$\alpha(\bar{H}/(Z_n(\bar{G})\cap\bar{H}))=\widetilde{H}/(Z_n(\widetilde{G})\cap
\widetilde{H}).$ Also one can see that
$\beta:[_n\bar{H},\bar{G}]\rightarrow[_n\widetilde{H},\widetilde{G}]$
by the rule $\beta(\bar{x})=\widetilde{x}$ is an isomorphism. By
Lemma 3.1, $(\alpha,\beta)$ is a relative $n$-isoclinism from
$(\bar{H},\bar{G})$ to $(\widetilde{H},\widetilde{G})$. $\square$\\

For proving (i) Theorem D, we need of the
following lemma.\\

{\underline{\bf{Lemma 3.4.}}}~ {\em Let $G$ be a group and $H$ be a
subgroup of $G$ such that $G=HZ_n(G)$. Then $(H,H) \
_{\widetilde{n}} \ (H,G) \ _{\widetilde{n}} \ (G,G)$.}\\

\indent {\it \bf Proof.}~ We want to prove $(H,H) \ _{\widetilde{n}}
\ (H,G)$. Let $G=HZ_n(G)$. We may easily see that $Z_n(H)=Z_n(G)\cap
H$. Thus  $H/Z_n(H)= H/(Z_n(G)\cap H)$ is isomorphic to
$HZ_n(G)/Z_n(G)=G/Z_n(G)$. Therefore $\alpha:H/Z_n(H)\rightarrow
G/Z_n(G)$ is an isomorphism which is induced by the inclusion
$i:H\rightarrow G.$ Furthermore, we can consider $\alpha$ as
isomorphism from $H/Z_n(H)$ to $H/Z_n(G)\cap H$.

On the other hand, $[_nH,G]=[_nH,HZ_n(G)]=\gamma_{n+1}(H)$. By Lemma
3.1, the pair $(\alpha,1_{\gamma_{n+1}(H)})$ allows us to state that
$(H,H) \ _{\widetilde{n}} \ (H,G)$.

The remaining cases $(H,H) \ _{\widetilde{n}} \ (G,G)$ and $(H,H) \
_{\widetilde{n}} \ (H,G)$ follow by a similar argument.$\square$\\

{\underline{\bf{Proof of Theorem D.}}}~ (i). By Lemma 3.4, we have
$(H,H) \ _{\widetilde{n}} \ (H,G) \ _{\widetilde{n}} \ (G,G)$.
$$d^{(n)}(H)=d^{(n)}(H,H)=d^{(n)}(H,G)=d^{(n)}(G,G)=d^{(n)}(G),$$ by
Theorem B and the result follows.

(ii). Assume $\varphi\in Aut(G)$. Then $\varphi$ induces the
isomorphisms $\alpha$ from $G/Z_n(G)$ to $G/Z_n(G)$ by the rule
$\alpha(gZ_n(G))=\varphi(g)Z_n(G)$ and $\beta$ from $[_nH,G]$ to
$[_n\varphi(H),G]$ by the rule
$\beta([h_1,...,h_n,x])=\varphi([h_1,...,h_n,x])$. Note that
\[\alpha(\frac{H}{Z_n(G)\cap H})=\frac{\varphi(H)}{Z_n(G)\cap\varphi(H)}.\] On another
hand, for every $g\in G$ and $h_i\in H$, $1\leq i\leq n$, we have $\varphi(g)\in\alpha(gZ_n(G)),$
$\varphi(h_i)\in\alpha(h_i(Z_n(G)\cap H))$ and
$\beta([h_1,...,h_n,g])=[\varphi(h_1),...,\varphi(h_n),\varphi(g)].$ By Lemma 3.1, the pair $(\alpha,\beta)$
implies that $(H,G) \ _{\widetilde{n}} \ (\varphi (H),G)$ and so  \[d^{(n)}(H,G)=d^{(n)}(\varphi(H),G).\ \ \ \
\square\]

Theorem C has two useful consequences, as we see in the next statements.\\

{\underline{\bf{Corollary 3.5.}}}~ {\em Let $H$ be subgroup of a
group $G$. Then there exists a group $G_1$ and a normal subgroup
$H_1$ of $G_1$ such that $(H,G) \ _{\widetilde{1}} \ (H_1,G_1)$ and
$Z(G_1)\cap H_1\subseteq H_1 \cap G'_1$. In particular, if $G$ is
finite, then
 $G_1$ is finite}.\\

\indent {\it \bf Proof.}~ Let $1\rightarrow R\rightarrow
F\rightarrow G\rightarrow1$ be a free presentation of $G$, $S$ be a
 subgroup of $F$, $H$ be a group isomorphic to $S/R$. If
$\bar{F}=F/(R\cap \bar{F}')$ and $\bar{S}=S/(R\cap \bar{F}')$, then
Theorem C with $n=1$ implies $(H,G) \ _{\widetilde{1}} \
(\bar{S},\bar{F})$. On another hand,
$(Z(\bar{F})\cap\bar{S})/(Z(\bar{F})\cap\bar{S} \cap \bar{F}')$ is
isomorphic to $((Z(\bar{F})\cap\bar{S})\bar{F}')/\bar{F}'$, which is
a subgroup of $ \bar{F}/\bar{F}'$. Therefore, for some normal
subgroup $\bar{B}$ of $\bar{F}$,
$Z(\bar{F})\cap\bar{S}=(Z(\bar{F})\cap \bar{S} \cap \bar{F}')\times
\bar{B}$.  Now $\bar{B}\cap \bar{F}'=1$ and we have $(H,G) \
_{\widetilde{1}} \ (H_1,G_1)$ again by Theorem C with $n=1$, where
$G_1=\bar{F}/\bar{B}$ and $H_1=\bar{S}/\bar{B}$. Furthermore
$Z(G_1)\cap H_1$ is isomorphic to $Z({\bar{F}/\bar{B})\cap
\bar{S}/\bar{B}=(Z(\bar{F})\cap\bar{S}})/\bar{B}$, which is a
subgroup of $(\bar{S} \cap \bar{F}')\bar{B}/\bar{B}=H_1 \cap G'_1$.
In particular, if $G$ is finite, then the index $|G_1:Z(G_1)\cap
H_1|$ is finite. But, $Z(G_1)\cap H_1\subseteq H_1 \cap G'_1,$
therefore $|G:G'_1|$ is finite and so $G_1$. $\square$ \\

{\underline{\bf{Corollary 3.6.}}}~ {\em Assume that $H$ is a
subgroup of a finite group $G$. Then there exists a finite group
$G_1$ and a normal subgroup $H_1$ of $G_1$ such that
$d(H,G)=d(H_1,G_1)$ and $Z(G_1)\cap H_1 \subseteq G'_1 \cap H_1.$}\\

\indent {\it \bf Proof.}~ By Theorem B and Corollary 3.5 the result
follows. $\square$\\

We know that $$D_8=\langle a,b \ | \ a^4=b^2=(ab)^2=1\rangle.$$ It
is easy to check that
$$(D_8,\langle a\rangle)\ _{\widetilde{_1}} \ (D_8,\langle
a^2,b\rangle) \ _{\widetilde{_1}} \ (D_8,\langle a^2,ab\rangle).$$
and that
$$d(D_8,\langle a\rangle)=d(D_8,\langle a^2,b\rangle)=d(D_8,\langle
a^2,ab\rangle)=\f34.$$ Theorem E shows that all pairs of groups with
the relative commutativity degree $\frac{3}{4}$ belong to the class
of relative 1-isoclinism of $( \langle a\rangle,D_8)$.\\

The following lemma gives an upper bound for $d(H,G)$ which will be
used in the proof of Theorem E.\\

{\underline{\bf{Lemma 3.5.}}}~ {\em For every subgroup $H$ of a
finite group $G$, \[d(H,G)\leq \f12(1+\f{|Z(G)\cup Z(H)|}{|G|}).\]}
\indent {\it \bf Proof.}~
\[\begin{array}{lcl}
d(H,G)&=&\ds\f1{|G||H|}|\{(h,g)\in H\times
G:[h,g]=1\}|=\f1{|G|}\sum_{g\in G}\f{|C_H(g)|}{|H|}\vspace{.2cm}\\
&=&\ds\f1{|G|}\left(\sum_{g\in Z(G)\cup
Z(H)}\f{|C_H(g)|}{|H|}+\sum_{g\notin Z(G)\cup
Z(H)}\f{|C_H(g)|}{|H|}\right)\vspace{.2cm}
\end{array}\]
\[\begin{array}{lcl}
&\leq&\ds\f1{|G|}\left(|Z(G)\cup Z(H)|+\f12(|G|-|Z(G)\cup Z(H)|)\right)\vspace{.2cm}\\
&=&\ds\f12(1+\f{|Z(G)\cup Z(H)|}{|G|}).\hspace{0.3cm} \square
\end{array}\]\\

{\textit{\bf{Proof of Theorem E.}}}~ Assume $d(H,G)=\f34$. Then $H$
is abelian by [4, Theorems 2.2 and 3.3] and $|G : H| \leq 2$ by
Lemma 3.7. Moreover, $|H/Z(G)|=2$ by [4, Theorem 3.10] and so
$|G:Z(G)|=4$. Therefore $G/Z(G)$ is a $2$-elementary abelian group
of rank 2 so we may define the isomorphism $\alpha$ from $G/Z(G)$ to
$D_8/Z(D_8)$ by $\alpha(\bar{x})=\bar{a}$ and
$\alpha(\bar{y})=\bar{b}.$ Since $Z(G)\subseteq H$, $H/Z(G)$ is
either $\langle \bar{x} \rangle$ or $\langle \bar{y} \rangle$ or
$\langle \bar{x}\bar{y} \rangle$.

Assume that $H/Z(G)=\langle \bar{x} \rangle$. Then
$\alpha(H/Z(G))=\langle a\rangle/ \langle a^2\rangle$ and
$[H,G]=\langle x,y\rangle$. Therefore $\beta:[H,G]\rightarrow
\langle a^2\rangle$ by $\beta([x,y])=[a,b]$ is an isomorphism. Hence
$(\alpha,\beta)$ is a relative isoclinism from $(H,G)$ to $(\langle
a\rangle,D_8)$ by Lemma 3.1. Now we have the remaining cases
$H/Z(G)=\langle\bar{y}\rangle$ and $H/Z(G)=\langle
\bar{x}\bar{y}\rangle$. If  $H/Z(G)=\langle\bar{y}\rangle$, then a
similar argument shows that $(H,G) \ _{\widetilde{1}} \ (\langle
a^2,b,D_8\rangle)$. If $H/Z(G)=\langle \bar{x}\bar{y}\rangle$, then
a similar argument shows that $(H,G) \ _{\widetilde{1}} \ (\langle
a^2,ab\rangle,D_8)$. There are no other cases so we deduce that $(H,
G) \ _{\widetilde{1}} \ (\langle a\rangle,D_8)$, as claimed.

Conversely, if $(H,G) \ _{\widetilde{1}} \ (\langle a\rangle,D_8)$,
then $d(H,G)=d(\langle a\rangle,D_8)=\f34$ and the result follows
from [4, Theorem 3.10].$\square$\\

Finally, we state the following conjecture. \\

{\underline{\bf{Conjecture.}}}~ \textit{Theorems B and
D hold when $G_1$ and $G_2$ are two infinite groups.}\\

We note that Definitions 2.1, 2.2 and 2.3 hold also in the infinite case. The same happens for Lemmas 3.1, 3.2,
3.4, Proposition 3.3, Theorem A and Theorem C. This allows us to ask whether conditions as in Theorems B and D
can happen in the infinite case. We strongly believe that the conjecture is true at least in the case of compact
groups, when a suitable notion of commutativity degree is introduced. We have some evidences in special cases as
[16, Theorems A and B] and  in recent submitted papers by the authors.

\section*{References}

\begin{itemize}
\item[]{[1]} \  J. C. Bioch, On n-isoclinic groups, \textit{Indag. Math.} {\bf{38}}
(1976), 400--407.
\item[]{[2]} \  J. C.  Bioch  and  R. W. van der Waall, Monomiality and
isoclinism of groups, \textit{J. Reine Ang. Math.} {\bf{298}}
(1978), 74--88.
\item[]{[3]} \ K. Chiti,   M. R. R. Moghaddam and A. R. Salemkar,
$n$-isoclinism classes and $n$-nilpotency degree of finite groups,
\textit{Algebra Colloq.} (2) \textbf{12} (2005), 225--261.
\item[]{[4]} \  A. Erfanian, P. Lescot  and R. Rezaei, On the relative commutativity
degree of a subgroup of a finite group, \textit{Comm. Algebra}
{\bf{35}} (2007), 4183--4197.
\item[]{[5]} \ P. X. Gallagher, The number of Conjugacy classes in
a finite group, \textit{Math. Z.} {\bf{118}} (1970), 175--179.
\item[]{[6]} \ W. H. Gustafson,  What is the probability that two groups
elements commute? ,  \textit{Amer. Math. Monthly} {\bf{80}} (1973),
1031--1304.
\item[]{[7]} \ M. Hall and J. K. Senior, \textit{The Groups of Order $2^n \ (n\leq
6)$}, Macmillan, New York, 1964.
\item[]{[8]} \ P. Hall, The classification of prime-power groups, \textit{J. Reine Ang. Math.} {\bf{182}} (1940),
130--141.
\item[]{[9]} \ P. Hall, Verbal and marginal subgroups, \textit{J. Reine Ang. Math.} {\bf{182}} (1940),
156--157.
\item[]{[10]} \ N. S. Hekster, On the structure of n-isoclinism classes
of groups, \textit{J. Pure Appl. Algebra} {\bf{40}} (1986), 63--85.
\item[]{[11]}  \ R. James, The groups of order $p^6$ ($p$ an odd
prime), \textit{Math. Comp.} {\bf{34}} (1980), 613--637.
\item[]{[12]} \ P. Lescot, Sur certains groupes finis, \textit{Rev. Math.
Sp\'eciales} {\bf{8}} (1987), 276--277.
\item[]{[13]} \ P. Lescot, Isoclinism classes and Commutativity
degrees of finite groups, \textit{J. Algebra} {\bf{177}} (1995),
847--869.
\item[]{[14]}  \ G. A.Miller, Relative number of non-invariant
operators in a group, \textit{Proc. Nat. Acad. Sci. USA} (2)
{\bf{30}} (1944), 25--28.
\item[]{[15]}  \  D. J. S. Robinson, \textit{A Course in the Theory of Groups}, Springer, Heidelberg, 1982.
\item[]{[16]}  \ F.G. Russo, A Probabilistic Meaning of Certain Quasinormal Subgroups, \textit{Int. J. Algebra} (1) \textbf{8} (2007),
385--392.
\end{itemize}
\end{document}